\documentclass[letterpaper,11pt]{article}

\usepackage{amsthm}
\usepackage{amsmath}
\allowdisplaybreaks[3]

\newtheorem{thm}{Theorem}[section]

\newtheorem{conj}[thm]{Conjecture}

\newcommand{\beq}[1]{\begin{equation}\label{#1}}
\newcommand{\enq}[0]{\end{equation}}

\newcommand{\mn}[0]{\medskip\noindent}
\newcommand{\nin}[0]{\noindent}

\newcommand{\sm}[0]{\setminus}

\newcommand{\m}[0]{{\cal M}}
\newcommand{\M}[0]{{\cal M}}

\newcommand{\T}[0]{{\cal T}}

\newcommand{\ra}[0]{\rightarrow}

\newcommand{\ZZ}[0]{{\bf Z}}
\newcommand{\fctu}[0]{{U}}

\newcommand{\ff}[0]{{\bf f}}

\newcommand{\sss}[0]{{\bf s}}

\newcommand{\YY}[0]{{\bf Y}}

\renewcommand{\qed}[0]{\begin{flushright} \rule{2mm}{3mm} \end{flushright}}

\def\qqqed{\null\nobreak\hfill\hbox{\rule{2mm}{3mm} }\par\smallskip}

\newcommand{\ga}[0]{\alpha }

\newcommand{\gd}[0]{\delta }
\newcommand{\gD}[0]{\Delta }

\newcommand{\gs}[0]{\sigma}

\usepackage[usenames]{color}

\newcommand{\sugg}[1]{}

\begin{document}

\renewcommand{\thefootnote}{\fnsymbol{footnote}}
\title{Asymptotics of the Upper Matching Conjecture\footnotemark }

\author{
L. Ilinca and J. Kahn}
\date{}

\footnotetext{ * Supported by NSF grant DMS0701175.}

\maketitle

\begin{abstract}
We give upper bounds for the number $\Phi_\ell(G)$ of matchings of 
size $\ell$ in
(i)
bipartite graphs $G=(X\cup Y, E)$ with specified degrees $d_x$
($x\in X$),
and
(ii) general graphs $G=(V,E)$ with all degrees specified.
In particular, for $d$-regular, $N$-vertex graphs, our bound is best
possible up to an error factor of the form
$\exp[o_d(1)N]$, where $o_d(1)\ra 0$ as $d\ra \infty$.
This represents the best progress to date on the ``Upper Matching Conjecture" of Friedland, Krop, Lundow and Markstr\"om.

Some further possibilities are also suggested.
\end{abstract}

\section{Introduction}

Throughout this paper $G$ will be a graph
without isolated vertices
on a vertex set $V$ of size $N$,
and $\ell$ an integer with $0 \leq \ell \leq N/2$.
When $G$ is bipartite the bipartition will be $X\cup Y$.
We write $\M_\ell(G)$ for the set of matchings of size $\ell$
(or $\ell$-{\em matchings}) of $G$, and set $\Phi_\ell(G)=|\M_\ell(G)|$.
The following ``Upper Matching Conjecture" was proposed by
Friedland {\em et al.}
 \cite{FKLM}.
\begin{conj}
\label{UMC}
If $G$ is $d$-regular and $2d \big| N$ then
\[
\Phi_\ell(G) \leq \Phi_\ell \left( \tfrac{N}{2d}K_{d,d} \right).
\]
\end{conj}

\nin
(As usual $tK_{d,d}$ is the union of $t$ disjoint
copies of the complete bipartite graph $K_{d,d}$.
A generalization of Conjecture \ref{UMC} that doesn't require
$2d\big| N$ was suggested by
Friedland, Krop and Markstr\"om \cite[(6.2)]{FKM}.)

\medskip
When $N$ is even and $\ell=N/2$, $\Phi_{N/2}(G)$ is the number,
$\Phi(G)$, of perfect matchings of $G$.
Here the bipartite case of Conjecture \ref{UMC} is
contained in the well-known theorem of
Br\'{e}gman \cite{Bre} (formerly the {\em Minc Conjecture} \cite{Minc}),
which in graph-theoretic language says that for $G$ bipartite on $X\cup Y$
(with $|X|=|Y|$)
\beq{Bregman}
\Phi(G) \leq \prod_{x\in X}  \left(d_x!\right)^{1/d_x}
\enq
(where $d_x$ is the degree of $x$). An analogous bound for general graphs was first observed in
\cite{KL} and rediscovered and/or reproved in \cite{AF, CR, Ego, Fri}.

Friedland {\em et al.} \cite{FKM} proved
Conjecture \ref{UMC} for $\ell=2$. Here we are really thinking of
$\ell=\Theta(N)$ and set $\ga=2\ell/N$,
the fraction of vertices covered by an $\ell$-matching. Carroll,
Galvin and Tetali \cite{CGT} provided some evidence in favor of
Conjecture \ref{UMC}, showing that (for $G$ as in the conjecture)
\beq{CGTub}
\log \Phi_\ell(G) \leq \frac{N}{2}
\big[
\ga \log d + H(\ga)
\big].
\enq
(Throughout this paper $\log$ is $\log_2$ and $H$ is the usual binary entropy function; for entropy basics see e.g. \cite{McE}.)
This contrasts with the
lower bound (which Conjecture \ref{UMC} would say is the truth) given by
\beq{UMClb}
\log \Phi_\ell(\tfrac{N}{2d}K_{d,d}) = \frac{N}{2}
\left[
\ga \log d + 2 H(\ga) + \ga \log \frac{\ga}{e} + o_d(1)
\right]
\enq
(where $o_d(1)\ra 0$ as $d\ra \infty$).

\medskip
Here we are primarily interested in closing the gap between
the parts of \eqref{CGTub} and \eqref{UMClb} that are on the order of $N$ as $d\rightarrow \infty$:

\begin{thm}
\label{dregularUB}
If $G$ is  $d$-regular, then
\[
\log \Phi_\ell(G) \leq \frac{N}{2}
\left[
\ga \log d + 2 H(\ga) + \ga \log  \frac{\ga}{e} +
\frac{\log d}{d-1}
\right].
\] 
\end{thm}
\nin
(Here and elsewhere we interpret $\log x/(x-1)$ as
$\log e$ when $x=1$.
Of course for Theorem \ref{dregularUB} we could simply
disallow the uninteresting case $d=1$, but the convention will
be helpful later.)

Theorem \ref{dregularUB}
is a special case of our main result, Theorem \ref{lMub}, which bounds $\Phi_\ell(G)$ in terms of the degree sequence $\{d_v\}_{v\in V}$ or, for $G$ bipartite on $X \cup Y$, $\{d_x\}_{x\in X}$.
For $W\subseteq V$, let
$d_\ell(W)=\sum\{\prod_{v\in S}d_v:S \in \binom{W}{\ell}\} $.
\begin{thm}
\label{lMub}
{\rm (i)} Suppose $G$ is bipartite (on $X\cup Y$)
and $0 \leq \ell \leq \min\{|X|, |Y|\}$, and set
$\gd_X =\min\{d_x: x \in X\}$ and $\ga_{_Y}=\ell/|Y|$. Then
\[
\log \Phi_\ell(G) \leq \log d_\ell(X)+{|Y|}
\left[
H(\ga_{_Y}) + \ga_{_Y} \log \frac{\ga_{_Y}}{e} + \frac{\log \gd_X}{\gd_{X}-1}
\right].
\]
{\rm (ii)} For a general $G$ (on V) with minimum degree $\gd$,
\beq{general}
\log \Phi_\ell(G) \leq  \frac{1}{2} \log d_{2\ell}(V) + \frac{N}{2}
\left[
H(\ga) + \ga \log \frac{\ga}{e} + \frac{\log \gd}{\gd-1}
\right].
\enq
\end{thm}
\nin
See the remark following \eqref{preremk}
for an explanation of the bound in (i).
Note that for $d$-regular $G$,
\[
\frac{1}{2} \log d_{2\ell}(V) = \frac{1}{2} \left( \log\binom{N}{2\ell}+2\ell \log d\right)<\frac{N}{2}\left( H(\ga)+ \ga \log d\right) ,
\]
so (ii) includes Theorem \ref{dregularUB}.

\medskip
The proofs of these results are mostly based on entropy considerations, in the spirit of Radhakrishnan's proof of Br\'{e}gman's Theorem \cite{Rad}, and, for example, the more recent \cite{CK}. Here, as for some related problems, most of the work deals with the bipartite case.
The passage to general graphs is then accomplished {\em via}
an easy correspondence between ordered pairs of $\ell$-matchings of $G$ and a subset of the $(2\ell)$-matchings of the bipartite double cover of $G$;
this correspondence, which goes back at least to Gibson \cite{Gib}, was
recently rediscovered by Alon and Friedland \cite{AF}.

The next section is devoted to the proof of Theorem \ref{lMub}.
In Section \ref{Remarks} we propose an extension of Br\'egman's bound
\eqref{Bregman} to unbalanced graphs, which would also be
a precise version of our main inequality \eqref{GSbound}.


\section{Proof of Theorem \ref{lMub}}

As noted above, (ii) will follow easily from (i).
We begin with the latter.
%
The main point here is establishing the bound when
$|X|=\ell$; that is, for $G$ bipartite on $S\cup Y$ with
$|S|=\ell$, $|Y|=M$ and $\gd_S=\min_{x\in S}d_x$,
\beq{GSbound}
\log \Phi_\ell(G) \leq \sum_{x\in S} \log d_x+M \left[H(\ga_{_Y})
+\ga_{_Y} \log \frac{\ga_{_Y}}{e}+ \frac{\log \gd_S}{\gd_S-1} \right].
\enq
This easily gives (i):
setting $G_S=G[S\cup Y]$ for $S\subseteq X$, we have
\begin{align}
\Phi_\ell(G) &= \sum_{S \in \binom{X}{\ell}} \Phi_\ell(G_S)
\nonumber\\
&\leq
 \sum_{S \in \binom{X}{\ell}} \Big( \prod_{x\in S}d_x \Big)  \exp_2
 \left[ M \left(H(\ga_{_Y}) +\ga_Y \log \frac{\ga_Y}{e}+
 \frac{\log \gd_S}{\gd_S-1} \right) \right]
\nonumber\\
&\leq
d_\ell(X) \exp_2 \left[ |Y| \left(H(\ga_{_Y})
+\ga_{_Y} \log \frac{\ga_{_Y}}{e}+ \frac{\log \gd_X}{\gd_{X}-1}
\right) \right].
\label{preremk}
\end{align}
\qed

\mn {\em Remark.}
The bound in \eqref{GSbound} (apart from the error term) is quite natural: $\exp_2 \left[MH(\ga_{_Y})\right]$ is roughly the number of ways to choose an $\ell$-subset $T$ of $Y$ to be used in the matching; on the other hand, for $T$ uniform from $\binom{Y}{\ell}$,
the ``typical" $T$-degree of $x\in S$ is $\ga_{_Y} d_x$, and
\[
\sum_{x\in S} \log d_x + M \ga_{_Y} \log \frac{\ga_{_Y}}{e} = \sum_{x\in S} \log \frac{\ga_{_Y} d_x}{e}
\]
is essentially Br\'egman's bound \eqref{Bregman} for these degrees.

\bigskip
 For the proof of \eqref{GSbound} we think of $\ell$-matchings of $G$
 as (injective) functions $f:S \rightarrow Y$, using
$R(f)$ for the range
of $f$ and $f_W$ for the restriction of $f$ to $W \subseteq S$.
For a permutation $\gs$ of $S$ (thought of as an ordering of $S$)
and $x\in S$, set $B(\gs, x)=\{w\in S: \gs(w) < \gs(x) \}$.
In what follows $x$ and $y$ range over $X$ and $Y$ respectively.
Expressions of the form $0\cdot\log b$ are always interpreted as zero.

Let $\ff $ be a random (uniform) $\ell$-matching of $G$ and
$p(x,y)= \Pr(\ff(x) = y)$.
Let $\sss$ be a random (uniform) permutation of $S$ and
$\YY_x= \ff_{B(\sss, x)}$. Thus, if we think of choosing $\ff$-values in
the order given by $\sss$, $\YY_x$ tells us what has happened prior to the
choice of $\ff(x)$. Our argument through \eqref{fZ} closely follows that
of \cite{CK}.
By the ``chain rule" for entropy, we have
\begin{align}
\log \Phi_\ell(G) &= H(\ff)
\nonumber\\
&= \frac{1}{\ell!} \sum_{\sigma} \sum_{x} H(\ff(x)|\ff_{B(\gs,x)})
\nonumber\\
&= \sum_{x} \sum_{\gs} \sum_{g} \frac{1}{\ell!} \Pr(\ff_{B(\gs,x)}=g)H(\ff(x)|\ff_{B(\gs,x)}=g)
\nonumber\\
&=
\sum_x H(\ff(x)|\YY_x),
\label{fY}
\end{align}
where $\gs$ ranges over the possible values of $\sss$ and $g$ over the possible values of $\ff_{B(\gs, x)}$.

Let $\ZZ_x=Y \, \sm \, \ff(B(\sss, x))$, the set of vertices of $Y$ that
remain available when we come to specify $\ff(x)$. Since $\ZZ_x$ is a
function of $\YY_x$, we have $H(\ff(x)|\YY_x) \leq H(\ff(x)|\ZZ_x)$ and
so, by \eqref{fY},
\beq{fYZ}
\log \Phi_\ell(G) \leq \sum_x H(\ff(x)| \ZZ_x).
\enq
\nin
In what follows we use $y$ and $Z$ for possible values
of $\ff(x)$ and $\ZZ_x$ for the $x$ under discussion,
in particular restricting to $y$'s
for which $p(x,y)\neq 0$.

From this point through \eqref{Fv1}, with the exception of
\eqref{Hf}, we fix $x\in S$.
Let $\Pr(Z)=\Pr(\ZZ_x=Z)$, $\Pr(Z|y)=\Pr(\ZZ_x=Z|\ff(x)=y)$ and so on,
and set $q_k(y) = \Pr(|\ZZ_x| = k | \ff(x) = y)$ and
$r_k(y) = \Pr (|\ZZ_x|=k, y\in \ZZ_x)$. Then, with
\(
F(x, y)= \sum_k q_k(y) \log \frac{r_k(y)}{q_k(y)}
\), we have

\begin{align}
H(\ff(x)|\ZZ_x) &= \sum_{Z}\Pr(Z)\sum_{y}\Pr(y|Z)\log \frac{1}{\Pr(y|Z)}
\nonumber\\
&= \sum_y \sum_Z \Pr(y,Z) \log \frac{\Pr(Z)}{\Pr(y,Z)}
\nonumber\\
&= \sum_y p(x,y) \left[ \log \frac{1}{p(x,y)} + \sum_Z \Pr(Z|y) \log \frac{\Pr(Z)}{\Pr(Z|y)} \right]
\nonumber\\
&= H(\ff(x))+ \sum_y p(x,y) \sum_Z \Pr(Z|y) \log \frac{\Pr(Z)}{\Pr(Z|y)}
\nonumber\\
&= H(\ff(x))+ \sum_y p(x,y) \sum_k q_k(y) \sum_{|Z|=k}
\frac{\Pr(Z|y)}{q_k(y)} \log \frac{\Pr(Z)}{\Pr(Z|y)}
\nonumber \\
&\leq H(\ff(x)) + \sum_y p(x,y) \sum_k q_k(y) \log \Bigg[ \frac{1}{q_k(y)} \sum_{\substack{|Z|=k\ \\ y\in Z }} \Pr(Z)\Bigg]
\nonumber\\
&=  H(\ff(x))+ \sum_y p(x,y) F(x,y),
\label{fZ}
\end{align}
the inequality following from concavity of the logarithm.
Rewriting
\[
H(\ff(x))=\log d_x - \sum_y p(x,y) \log (d_x p(x,y))
\]
and applying \eqref{fYZ} and \eqref{fZ} (and momentarily
unfixing $x$) gives
\beq{Hf}
\log \Phi_\ell(G) \leq \sum_{x\in S} \log d_x + \sum_x \sum_y p(x,y) \big[ F(x,y) -  \log (d_x p(x,y)) \big].
\enq
The main part of the proof involves bounding the second term in \eqref{Hf},
and in particular $F(x,y)$.

(We again fix $x$.)
For $y\in Y$ set $\mu_y = \Pr(y \in R(\ff))$ and $\nu_y = 1-\mu_y$.
%
Since $|\ZZ_x|$ and $\ff$ are independent, we have
\beq{qk}
q_k(y) = \Pr(|\ZZ_x|=k)=
\left\{
\begin{array}{l l}
  0 & \quad \text{for } k \leq M-\ell ,\\
  1/\ell & \quad \text{for } M-\ell + 1\leq k \leq M,\\
\end{array} \right.
\enq
while $\Pr(y\in \ZZ_x \,|\, y\not \in R(\ff) \sm \{\ff(x)\})=1$ and  $\Pr(y\in \ZZ_x \,|\, y\in R(\ff) \sm \{\ff(x)\},|\ZZ_x|=k )= (k-(M-\ell)-1)/(\ell -1)$. Thus
\begin{align}
r_k(y) & = \Pr(|\ZZ_x|=k)
[
\Pr (y\in R(\ff) \sm \{\ff(x)\}) \cdot \tfrac{k-(M-\ell)-1}{\ell -1}
\nonumber\\
& ~~~~~~~~~~~~~~~~~~~~~~~~~~~~~~~~~~~~~~~~~~~~~~+ \Pr(y\not \in R(\ff) \sm \{\ff(x)\}) ]
\nonumber\\
& = q_k(y) \left[ \big( \mu_y - p(x,y) \big) \tfrac{k-(M-\ell)-1}{\ell -1} + \big( \nu_y + p(x,y) \big) \right]
\nonumber
\end{align}
and (for $M-\ell + 1\leq k \leq M$)
\beq{rk_over_qk}
\frac{r_k(y)}{q_k(y)}  =
\left( \mu_y - p(x,y) \right) \tfrac{k-(M-\ell)-1}{\ell -1} +
\left(\nu_y + p(x,y)\right).
\enq

Let
\[
\fctu(t)=\log \left[\left( \mu_y - p(x,y)\right) t + \left( \nu_y + p(x,y) \right) \right]    ~~~~~(t\in [0,1])
\]
and
\[
f(t)=
\left\{
\begin{array}{l l}
 \frac{t}{1-t} \log \frac{1}{t} \, & \quad \text{if } t \in (0,1) ,\\
0 & \quad \text{if } t=0,\\
  \log e & \quad \text{if } t=1.
\end{array}
\right.
\]

In view of \eqref{qk} and \eqref{rk_over_qk}, we have
(with justification of \eqref{Fv1} to follow),
\begin{align}
F(x,y) ~& \leq
\sum_{k=M-\ell+1}^{M} \frac{1}{\ell} \, \fctu \!\left( \tfrac{k-(M-\ell)-1}{\ell -1} \right)
~ = ~
\sum_{j=0}^{\ell -1} \frac{1}{\ell} \, \fctu \!\left( \tfrac{j}{\ell -1} \right)
\label{Fv0}\\
& \leq
~\int_0^1 \fctu(t)dt ~
=~
G(x,y) - \log e,
\label{Fv1}
\end{align}
where $G(x,y)=f(\nu_y+p(x,y))$.
The equality in line \eqref{Fv1} is trivial if
$\mu_y=p(x,y)$, and otherwise is given by
the fact that for $a\not=0$ and $b=1-a>0$,
\[
\int_0^1 \log(at+b) dt = \frac{\log e}{a}
\Big[
(at+b) \ln(at+b)-at \big|_{t=0}^{t=1}
\Big]
= \frac{b}{a} \log \frac{1}{b} - \log e.
\]
The inequality in \eqref{Fv1} requires only the concavity of $U$,
as follows.  Let $U^*$ be the smallest concave function that
agrees with $U$ at the points $j/(\ell-1)$; namely,
for $1\leq j\leq\ell-1$ and
$\tfrac{j-1}{\ell-1}\leq x\leq \tfrac{j}{\ell-1}$,
$$U^*(x) = (j-(\ell-1)x)U(\tfrac{j-1}{\ell-1})
+((\ell-1)x-(j-1))U(\tfrac{j}{\ell-1}) .
$$
Then $U^*\leq U$ and, setting
$a_i = U(i/(\ell-1))$, we have
\begin{eqnarray*}
\int_0^1U^*(t)dt &= &
\tfrac{1}{2(\ell-1)}[(a_0+a_1)+\cdots + (a_{\ell-2}+a_{\ell-1})]\\
&=&
\tfrac{1}{\ell-1}(a_0+\cdots +a_{\ell-1})
-\tfrac{1}{2(\ell-1)}(a_0+a_{\ell -1})\\
&\geq& \tfrac{1}{\ell}(a_0+\cdots a_{\ell-1}),
\end{eqnarray*}
which is the right hand side of \eqref{Fv0}.
(The inequality, which is equivalent to
$2(a_0+\cdots a_{\ell-1})\geq \ell(a_0+a_{\ell-1})$,
follows from the concavity of $U$
(an instance of Karamata's Inequality, e.g. \cite{HLP}).)\qqqed

\medskip
Thus, now letting $x$ vary, we have
\begin{align}
\sum_x \sum_y p(x,y)F(x,y)
& \leq
\sum_x
\sum_y
p(x,y)
\left[ G(x,y) - \log e
\right]
\nonumber\\
& =  \sum_x \sum_y p(x,y) G(x,y) - \ell \log e
\nonumber\\
&=
\sum_x \sum_y p(x,y) G(x,y) - M\ga_{_Y} \log e .
\label{Fv2}
\end{align}

\mn
We will approximate the last sum by
\begin{align}
\gD & :=
\sum_x \sum_y p(x,y) f(\nu_y)
=
\sum_y \nu_y \log \frac{1}{\nu_y}
\nonumber
\\
&
\leq -M (1-\ga_{_Y})\log (1-\ga_{_Y})
\nonumber
\\
& =
M [H(\ga_{_Y})+ \ga_{_Y} \log \ga_{_Y}] \label{middle}
\end{align}
\nin
(where we used $\sum_x p(x,y)=\mu_y$ in the first line
and Jensen's Inequality in the second).
Adding and subtracting $\gD$ from the right side of
\eqref{Fv2} and using \eqref{middle} gives
\begin{eqnarray}\label{Fv3}
\sum_x \sum_y p(x,y) F(x,y)
&\leq &
\sum_x \sum_y p(x,y)
[ G(x,y)
- (\nu_y/\mu_y) \log (1/\nu_y)]
\nonumber\\
&&~~~~~~~
+ M \left[H(\ga_{_Y})+ \ga_{_Y} \log (\ga_{_Y}/e)
 \right] .
\end{eqnarray}

\nin
We next observe that
$$
G(x,y)
- (\nu_y/\mu_y) \log (1/\nu_y) = f(\nu_y+p(x,y))-f(\nu_y) \leq
f(p(x,y)),
$$
the inequality holding because $f$ is concave with $f(0)=0$.
Combining this with \eqref{Hf} and \eqref{Fv3} gives
\begin{align}
\log \Phi_\ell(G) \leq {} & \sum_x \sum_y p(x,y) \Big[ f(p(x,y))-
\log (d_x p(x,y)) \Big]
\nonumber\\
& ~~ + \sum_x \log d_x +
M \left[H(\ga_{_Y})+ \ga_{_Y} \log \frac{\ga_{_Y}}{e}
\right].
\label{inequality}
\end{align}
Finally, for $x\in S$ and $t\in (0,1]$, set
\[
g_x(t)=
\left\{
\begin{array}{l l}
 \frac{t^2}{1-t}\log \frac{1}{t}-t \log(d_x t) & \quad \text{if } t \not = 1, \\
  \log e -\log d_x& \quad \text{if } t=1.\\
\end{array}
\right.
\]
The double sum in \eqref{inequality} is then
$\sum_x \sum_y g_x(p(x,y))$, and
\eqref{GSbound} follows from (\eqref{inequality} and)
the observation that (since $g_x(t)=f(t)-t\log d_x$
is concave in $t$)
\[
\sum_y g_x(p(x,y))
\leq
d_x g_x(1/d_x)=
\frac{\log d_x}{d_x-1}
\leq \frac{\log \gd_S}{\gd_S-1} .
\]
\qed

\medskip
Finally we turn to (ii).
We consider the bipartite double cover, say $K$,
of $G$; that is, the
graph on $V \times \{0,1\}$ with edge set $\{(x,0)(y,1): xy \in E(G)\}$.
This is a bipartite graph on $2N$ vertices
(with $d_K(x,0)=d_K(x,1)=d_G(x)$), so, as shown above,
\beq{bicover}
\log \Phi_{2\ell}(K) \leq \log d_{2\ell}(V) + N \left[H(\ga) +\ga \log \frac{\ga}{e}+ \frac{\log \gd}{\gd-1} \right],
\enq
where, again, $\gd=\min \{d_v: v\in V\}$.

Let $\T^*$ be the set of $(2\ell)$-edge multisubgraphs of $G$ (that is, multigraphs whose underlying graphs are subgraphs of $G$) whose components are paths and cycles (possibly of length $2$), and let $\T$ consist of those members of $\T^*$ that do not contain odd cycles. For $T \in \T^*$ let $c(T)$ be the number of components of $T$ that are not $2$-cycles.

The natural projection $\psi((x,i))=x$ maps $\M_{2\ell}(K)$ surjectively to $\T^*$, with $|\psi^{-1}(T)|=2^{c(T)}$ for all $T\in \T^*$. On the other hand, $\varphi: \M_\ell(G) \times \M_\ell(G) \rightarrow \T$ given by $\varphi(M, M')=M\cup M'$ ({\it multiset} union) is a surjection with $|\varphi^{-1}(T)|=2^{c(T)}$ for all $T\in \T$. Thus we have

\[
\Phi_\ell(G)^2 = \sum_{T \in \T} 2^{c(T)} \leq \sum_{T \in \T^*} 2^{c(T)}
= \Phi_{2\ell}(K),
 \]
which, combined with \eqref{bicover}, gives \eqref{general}.\qed

\section{A conjecture}\label{Remarks}
In closing
we would like to propose a precise version of \eqref{GSbound}
that was
one of the original reasons for our interest in the present material.

\medskip
For $d\geq t>0$ with $d$ an integer,
set $\psi (d,t) = t^{-1}[\log d!-\log \Gamma(d-t+1)]$
(with $\Gamma$ the usual gamma function).

\begin{conj}\label{genMinc}
Let
$G$ be bipartite on $X\cup Y$ with
$|X|=\ell$ and $|Y|=M\geq \ell$, and for $x\in X$,
let $t_x=\ell d_x/M$.
Then
\beq{wild1}
\log\Phi_\ell(G)\leq \sum_{x\in X}\psi(d_x,t_x).
\enq
\end{conj}

\nin
Notice that for $t\in \ZZ$, one has
$\psi(d,t) = t^{-1}\log (d)_t$
(where, as usual, $(d)_t = d(d-1)\cdots (d-t+1)$), and
in particular $\psi(d,d)=d^{-1}\log d!$.
Thus
Conjecture \ref{genMinc} for $M=\ell$ is Br\'egman's Theorem,
and the full conjecture is a natural generalization thereof
which is sharp whenever $G$ is a disjoint union of complete
bipartite graphs $K_{X_i,Y_i}$ (where
$X=\cup X_i$ and $Y=\cup Y_i$) with $|X_i|/|Y_i|=\ell/M$ $\forall i$.
(Of course one can also
think of this in the original setting of the Minc Conjecture,
viewing it as a bound on the ``generalized permanent" of a
not-necessarily-square
$\{0,1\}$-matrix with given row sums.)

\medskip
To be honest, we haven't thought much about plausibility of
Conjecture \ref{genMinc} in its full generality.
We do feel pretty sure that it is true, for example, when
$d_x=d$ for every $x$ and $t:=\ell d/M$ (the average of the $d_y$'s)
is an integer.  Curiously, we can prove this when
$d_y=t ~\forall y$, which ought to be the worst case, but so far
not in general.

For an even wilder possibility, set, for $r\geq 0$ and $0<t\leq 2^r$,
$\varphi (r,t) = t^{-1}[\log\Gamma(2^r+1)-\log(2^r -t+1)]$.
Could it be that for $G$ as in
Conjecture \ref{genMinc},
$\ff$ random (but {\em not necessarily uniform}) from $\m_\ell(G)$,
and $t_x= (\ell/M) 2^{H(\ff(x))}$ ($x\in X$),
one has
\beq{wild}
H(\ff)\leq \sum_{x\in X} \varphi (H(\ff(x)),t_x)?
\enq
Notice that when $\ff(x)$ is uniform from the neighbors of $x$,
$2^{H(\ff(x))}$ is just $d_x$,
so
\eqref{wild} strengthens \eqref{wild1}.
The case $\ell=M$ of \eqref{wild} was suggested in \cite{CK}.

\medskip
\begin{flushleft}
Department of Mathematics\\
Indiana University\\
Bloomington, IN 47401\\
ilinca@indiana.edu\\
\end{flushleft}

\begin{flushleft}
Department of Mathematics\\
Rutgers University\\
Piscataway NJ 08854 USA\\
jkahn@math.rutgers.edu
\end{flushleft}
\end{document}